\documentclass[conf]{new-aiaa}

\usepackage[utf8]{inputenc}

\usepackage{graphicx}
\usepackage{amsmath}
\usepackage[version=4]{mhchem}
\usepackage{siunitx}
\usepackage{longtable,tabularx}
\usepackage{xcolor}
\usepackage{booktabs}
\usepackage{algorithm}
\usepackage{algpseudocode}
\usepackage[caption=false,font=small]{subfig}

\graphicspath{{Figs/}}

\newtheorem{remark}{Remark}

\title{Nonlinear Guidance for Arrival Time and Arrival Angle Control Using Trajectory Shaping}

\author{
Kun Wang\footnote{Assistant Professor, Department of Aerospace Engineering. }
and Andr\'es Infante Adri\'an\footnote{Ph.D. Student, Department of Aerospace Engineering.}
}

\affil{Universidad Carlos III de Madrid, Legan\'es 28911, Spain}

\begin{document}

\maketitle

\begin{abstract}
This paper proposes a nonlinear trajectory shaping guidance strategy for arrival time and angle control of a constant-speed fixed-wing unmanned aerial vehicle. The look angle is parameterized by a fourth-order polynomial. The nonlinear guidance problem is transformed into solving two coupled nonlinear integral equations with respect to two unknown guidance parameters. Directly solving these equations is challenging due to the strong coupling between the parameters. To address this, a two-stage solution procedure is developed. In the first stage, an analytical warm start is constructed. Specifically, by applying approximations, a linear relationship between the two guidance parameters is obtained. With this relationship, a scalar quadratic equation in the first guidance parameter is derived, leading to a good initial guess for the first guidance parameter. In the second stage, this initial guess is refined by solving a one-dimensional nonlinear equation. The obtained solution is finally used as the initial guess for solving the original two-dimensional nonlinear system of equations. Numerical simulations demonstrate that for the test cases, the obtained solution is very close to the open-loop optimal solution, even in highly nonlinear scenarios where some existing methods fail to generate a feasible solution.
\end{abstract}

\section{Introduction}

Guidance laws play a critical role in autonomous missions, as they directly determine the vehicle's maneuvering commands to meet the mission objectives. In practice, in addition to guiding the vehicle to a desired location, 
vehicles are often required to satisfy additional complex constraints. For instance, in time-critical missions, the vehicle may be required to reach the destination at a specific time to synchronize with other vehicles \cite{shanmugavel2010co,chen2023elongation}. In some scenarios, the vehicle may also need to arrive at the destination with a specific direction to achieve superior performance \cite{shaferman2008linear}. Autonomous missions may also require simultaneous control of both arrival time and arrival angle, as this can satisfy more complex practical requirements \cite{hou2023optimal}. However, these fully constrained requirements and the nonlinear kinematics make the design of guidance laws particularly challenging. 
As a result, significant research efforts have been devoted to the design of guidance laws for arrival-time-and-arrival-angle control in recent years. Based on whether the guidance law design requires linearization of the highly nonlinear kinematics, existing methods can be broadly categorized into linear and nonlinear guidance laws.

Most early works were based on linearized kinematics \cite{lee2007guidance,kim2013augmented,zhang2013guidance}. For instance, Lee et al. \cite{lee2007guidance} introduced a bias term for the arrival-time error in a minimum-jerk guidance law to precisely adjust the arrival time and arrival angle. Zhang et al. \cite{zhang2013guidance} proposed a guidance law composed of a
constructed biased proportional navigation law with an arrival angle constraint and a feedback control term for the arrival-time error.
In addition, linear optimal control was also employed to derive analytical optimal guidance laws \cite{chen2019optimal}. 

For linear guidance laws, one essential concept is time-to-go estimation, which is needed to regulate the arrival time. However, the accuracy of such estimation is often compromised by the linearization, and the corresponding guidance law may fail to generate a feasible solution when the engagement geometry is highly nonlinear. To this end, more recent works have focused on nonlinear guidance law design \cite{wang2022nonlinear,majumder2023three,singh2025terminal}. These methods typically rely on sliding mode control \cite{kumar2018sliding}, nonlinear optimal control \cite{wu2025nonlinear}, trajectory shaping \cite{kang2019generalized}, or virtual-target-based methods \cite{hou2023optimal,hu2018new,zhang2022virtual}. It is worth mentioning that solving the corresponding nonlinear optimal control problem may suffer from convergence issues \cite{wang2024physics}; therefore, machine learning has been combined with optimal control to generate the optimal guidance law in real time \cite{wang2025nonlinear}. As for the trajectory shaping method, the idea is to parameterize the state or control trajectory using certain functions, such as polynomials, and then embed the boundary conditions. As a result, the guidance problem, without any linearization, is transformed into solving a set of nonlinear equations for the unknown shaping parameters. Consequently, there is no need for time-to-go estimation, and the resulting guidance law can be applied to highly nonlinear scenarios. Therefore, trajectory shaping has been applied to arrival time control \cite{tekin2017polynomial,wang2026FOV}, arrival angle control \cite{wang2026trajectory}, and simultaneous arrival time and angle control \cite{kang2019generalized}. However, it was shown in \cite{wang2026ITCG} that solving the resulting nonlinear equations may be nontrivial.

In light of the above discussion, we develop a nonlinear guidance law for simultaneous arrival time and angle control based on trajectory shaping. 
In this way, there is no need for time-to-go estimation, and the resulting guidance law can be applied to highly nonlinear scenarios. 
Compared with the existing work in \cite{kang2019generalized}, the main contribution of this paper is twofold. 
First, the scaling techniques in \cite{wang2022nonlinear,wang2026FOV} are applied, which simplifies the derivation of the guidance law. 
Second, a carefully designed two-stage solution procedure is developed to solve the resulting nonlinear equations, thereby improving convergence. 
The simulation results show that the proposed method can generate a solution that is very close to the open-loop optimal solution, even in highly nonlinear scenarios where some existing methods fail to generate a feasible solution.

The rest of the paper is organized as follows. Section~\ref{sec:problem_formulation} formulates the nonlinear guidance problem. Section~\ref{sec:guidance_law_design} presents the design of the trajectory shaping guidance law and the two-stage solution procedure. Section~\ref{sec:simulations} provides numerical simulations to validate the proposed method. Finally, Section~\ref{sec:conclusions} concludes the paper.
\section{Problem Formulation}\label{sec:problem_formulation}
\begin{figure}[htbp]
      \centering
      \includegraphics[width=0.50\textwidth]{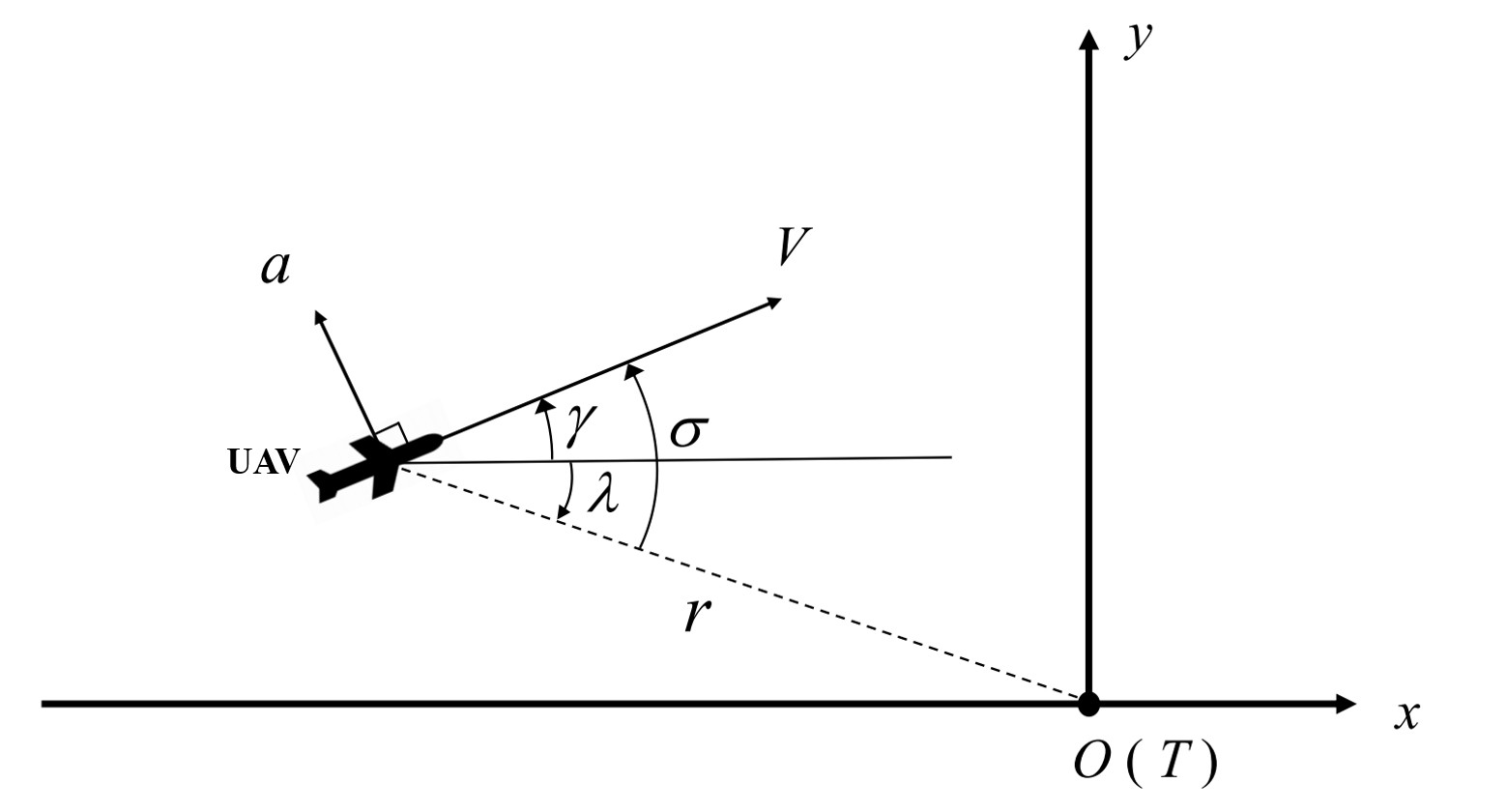}
      \caption{{Coordinate frame for the UAV and destination.}}
      \label{fig:engagement_geometry}
\end{figure}

Consider a two-dimensional engagement geometry where an unpowered fixed-wing UAV flies toward a stationary destination located at the origin. The UAV is modeled as a point mass moving at a constant speed \(V\). Without loss of generality, the speed is normalized as \(V=1\) \cite{wang2026FOV}. The engagement geometry is shown in Fig.~\ref{fig:engagement_geometry}.
Let \(r(t)\) denote the relative range, \(\lambda(t)\) the line-of-sight (LOS) angle, \(\gamma(t)\) the flight-path angle, and \(\sigma(t)\) the look angle. These angles satisfy
    ${\gamma(t)=\lambda(t)+\sigma(t).}$
{
The directions of \(\lambda(t)\), \(\sigma(t)\), and \(\gamma(t)\) are defined as positive when measured counterclockwise. With the normalized speed \(V=1\), the polar-coordinate dynamics are written as
}
\begin{equation}
{
\left\{
\begin{aligned}
    \dot r(t) &= -\cos\sigma(t),\\
    \dot\lambda(t) &= -{\sin\sigma(t)}/{r(t)},\\
    \dot\gamma(t) &= a(t),
\end{aligned}
\right.
}
\label{eq:polar_dynamics}
\end{equation}
{
where \(a(t)\) is the normalized normal acceleration command acting perpendicular to the UAV velocity vector. 
The state vector is therefore defined as
$
[r(t),\lambda(t),\gamma(t)]^T,
$
where \(r(t)\), \(\lambda(t)\), and \(\gamma(t)\) can be measured by a ranging sensor, a target seeker, and an inertial measurement unit, respectively. The initial condition at \(t=0\) is given by
}
\begin{equation}
{
    r(0)=r_0,\qquad
    \lambda(0)=\lambda_0,\qquad
    \gamma(0)=\gamma_0.
}
    \label{eq:initial}
\end{equation}
{
The objective is to guide the UAV to the destination at a prescribed arrival time \(t_f\), while achieving a desired arrival angle with bounded terminal acceleration. The terminal conditions are therefore specified as
}
\begin{equation}
{
    r(t_f)=0,\qquad
    \gamma(t_f)=\gamma_f,\qquad
    \sigma(t_f)=0,
}
    \label{eq:terminal}
\end{equation}
{
where \(\gamma_f\) is the desired arrival angle. The condition \(\sigma(t_f)=0\) ensures that the UAV velocity direction is aligned with the LOS direction at arrival. This guarantees a bounded terminal acceleration.
According to \cite{wang2026ITCG}, any desired arrival time \(t_f\) can be normalized to \(t_f=1\). The terminal conditions in \eqref{eq:terminal} are equivalently expressed as
}
\begin{equation}
{
    r(1)=0,\qquad
    \lambda(1)=\gamma_f,\qquad
    \sigma(1)=0.
}
    \label{eq:eq_terminal}
\end{equation}
{
Formally, the nonlinear guidance problem is to synthesize a control input \(a(t)\) that drives the UAV from \eqref{eq:initial} to \eqref{eq:eq_terminal}. Due to the nonlinear dynamics, this fully constrained guidance problem is computationally challenging for real-time onboard implementation. In the next section, a guidance law based on look angle shaping is developed to address this problem.
}
\begin{remark}
{
In practice, the normal acceleration command is subject to saturation limits due to the structural and aerodynamic constraints of the UAV. Explicitly incorporating these saturation limits into the analytical guidance-law design would significantly complicate the derivation. Therefore, they are not considered in the theoretical development, but are enforced in the numerical simulations.
}
\end{remark}

\section{Guidance Law Design} \label{sec:guidance_law_design}
\subsection{Parameterization of Look Angle Profile}
{
The look angle profile is prescribed as the following fourth-order polynomial:
}
\begin{equation}
    \sigma(t)=(t-1)^2\left(\sigma_0+\kappa_1t+\kappa_2t^2\right),
    \label{eq:shape_func}
\end{equation}
{
where \(\sigma_0=\sigma(0)\) is the initial look angle, and \(\kappa_1,\kappa_2\) are the guidance parameters to be determined. This structure inherently guarantees \(\sigma(1)=0\) and \(\dot{\sigma}(1)=0\). By applying L'H\^{o}pital's rule to \eqref{eq:polar_dynamics} as \(r\to0\), it follows that \(\dot{\lambda}(1)=0\). Hence, the terminal acceleration satisfies
\[
    a(1)=\dot{\gamma}(1)=\dot{\sigma}(1)+\dot{\lambda}(1)=0,
\]
which is desirable for smooth arrival.

The nonlinear guidance problem is thus reduced to determining \((\kappa_1,\kappa_2)\) such that the following nonlinear system of integral equations is satisfied:
\begin{equation}
{
\mathbf{F}(\kappa_1,\kappa_2)=
\begin{bmatrix}
r_0-\int_0^1\cos\sigma(t;\kappa_1,\kappa_2)\,dt\\
\lambda_0-\gamma_f-\int_0^1\dfrac{\sin\sigma(t;\kappa_1,\kappa_2)}{r(t)}\,dt
\end{bmatrix}
=\mathbf{0}.
}
\label{eq:exact_system}
\end{equation}
where
$
r(t;\kappa_1,\kappa_2)
=
r_0-\int_0^t \cos\sigma(s;\kappa_1,\kappa_2)\,ds .
$
{
Once the guidance parameters are obtained, the guidance command can be readily obtained. However, solving \eqref{eq:exact_system} directly is nontrivial because the two parameters are strongly coupled through nonlinear integrals. To improve convergence and computational efficiency, a two-stage procedure is developed below.
}

\subsection{Analytical Warm Start}
{
We first construct an analytical warm start by using some approximations. 
It should be emphasized that these approximations are used only to construct an initial guess. The final guidance parameters are still obtained by solving the exact nonlinear system of equations in \eqref{eq:exact_system}.
Under the small-angle and near-linear range-depletion approximations
}
\begin{equation}
{
    \sin\sigma\approx\sigma,\qquad r(t)\approx r_0(1-t),
}
\label{eq:approximations}
\end{equation}
{
the terminal LOS condition \(\lambda(1)=\gamma_f\) yields
}
\begin{equation}
{
    \lambda_0-\gamma_f
    =
    \frac{1}{r_0}\int_0^1
    (1-t)\left(\sigma_0+\kappa_1t+\kappa_2t^2\right)dt.
}
\end{equation}
{
Evaluating the integral gives
}
\begin{equation}
{
    \lambda_0-\gamma_f
    =
    \frac{1}{r_0}
    \left(
    \frac{1}{2}\sigma_0+
    \frac{1}{6}\kappa_1+
    \frac{1}{12}\kappa_2
    \right).
}
\end{equation}
Therefore, the two parameters satisfy the linear relationship
\begin{equation}
{
    \kappa_2=\Gamma-2\kappa_1,
}
    \label{eq:kappa2_linear}
\end{equation}
{
where \(\Gamma=12r_0(\lambda_0-\gamma_f)-6\sigma_0\).

Substituting \eqref{eq:kappa2_linear} into \eqref{eq:shape_func} gives the decoupled expression
}
\begin{equation}
{
    \sigma(t)=\sigma_0P_1(t)+\Gamma P_2(t)+\kappa_1P_3(t),
}
    \label{eq:sigma_decoupled}
\end{equation}
{
where \(P_1(t)=(t-1)^2\), \(P_2(t)=t^2(t-1)^2\), and \(P_3(t)=(t-1)^2(t-2t^2)\).  Applying the second-order approximation
}
\begin{equation}
{
    \cos\sigma\approx1-\frac{1}{2}\sigma^2
}
\label{eq:cos_taylor}
\end{equation}
{
to the range equation in \eqref{eq:exact_system} leads to the scalar quadratic equation
}
\begin{equation}
{
    A\kappa_1^2+B\kappa_1+C=0,
}
    \label{eq:quadratic_final}
\end{equation}
{
where
}
\begin{subequations}
\begin{align}
    A &= \int_0^1 P_3^2(t)\,dt
    = \frac{1}{630}, \\
    B &= 2\sigma_0\int_0^1P_1(t)P_3(t)\,dt
    +2\Gamma\int_0^1P_2(t)P_3(t)\,dt
    = \frac{\sigma_0}{35}+\frac{\Gamma}{1260}, \\
    C &= \sigma_0^2\int_0^1P_1^2(t)\,dt
    +2\sigma_0\Gamma\int_0^1P_1(t)P_2(t)\,dt
    +\Gamma^2\int_0^1P_2^2(t)\,dt
    -2(1-r_0) \nonumber\\
    &=
    \frac{\sigma_0^2}{5}
    +\frac{2\sigma_0\Gamma}{105}
    +\frac{\Gamma^2}{630}
    -2(1-r_0).
    \label{eq:coeff_C}
\end{align}
\end{subequations}

Let \(\Delta=B^2-4AC\) denote the discriminant of \eqref{eq:quadratic_final}. If \(\Delta \le 0\), the quadratic approximation does not provide two real roots. In this case, the vertex of the quadratic approximation, \(-B/(2A)\), is used as the initial estimate.
If \(\Delta>0\), two real candidates $$\kappa_1^{(1)}=\frac{-B-\sqrt{\Delta}}{2A}\quad\text{and}\quad\kappa_1^{(2)}=\frac{-B+\sqrt{\Delta}}{2A}$$ are obtained, and the candidate with the smaller approximate control effort is selected. Specifically, using the approximations in \eqref{eq:approximations}, the approximate normal acceleration used for selecting the initial guess is
\begin{equation}
{
    a(t)\approx \hat a(t)
    =
    \dot{\sigma}(t)-\frac{\sigma(t)}{r_0(1-t)}.
}
\label{eq:accel_linear}
\end{equation}
{
Using \eqref{eq:sigma_decoupled}, it can be written as
}
\begin{equation}
{
    \hat a(t;\kappa_1)
    =
    \sigma_0\tilde P_1(t)+
    \Gamma\tilde P_2(t)+
    \kappa_1\tilde P_3(t),
}
\label{eq:accel_decoupled}
\end{equation}
where
\begin{subequations}
\begin{align}
    \tilde P_1(t)&=\dot P_1(t)-\frac{P_1(t)}{r_0(1-t)},\\
    \tilde P_2(t)&=\dot P_2(t)-\frac{P_2(t)}{r_0(1-t)},\\
    \tilde P_3(t)&=\dot P_3(t)-\frac{P_3(t)}{r_0(1-t)}.
\end{align}
\end{subequations}
{
As \(P_1(t)\), \(P_2(t)\), and \(P_3(t)\) all contain the factor \((t-1)^2\), the apparent singularity at \(t=1\) is canceled. 
The approximate control effort is then defined as
\begin{equation}
    \hat J(\kappa_1)=\int_0^1\hat a^2(t;\kappa_1)\,dt.
\label{eq:control_effort_J}
\end{equation}
Since \(\hat a(t;\kappa_1)\) is a polynomial in \(t\) and is affine with respect to \(\kappa_1\), the above integral can be evaluated analytically as
\begin{equation}
    \hat J(\kappa_1)
    =
    Q_2\kappa_1^2+Q_1\kappa_1+Q_0,
    \label{eq:Jhat_quadratic}
\end{equation}
where
\begin{subequations}
\begin{align}
    Q_2 &= \frac{16r_0^2+r_0+1}{210r_0^2},\\
    Q_1 &=
    -\frac{\sigma_0(r_0-1)(2r_0+1)}{30r_0^2}
    -\frac{\Gamma(2r_0^2+r_0+1)}{210r_0^2},\\
    Q_0 &=
    \frac{\sigma_0^2(2r_0+1)^2}{3r_0^2}
    -\frac{\sigma_0\Gamma(r_0-1)(2r_0+1)}{15r_0^2}
    +\frac{\Gamma^2(2r_0^2+r_0+1)}{105r_0^2}.
\end{align}
\end{subequations}

The Taylor initial guess is finally selected as
\begin{equation}
{
    \kappa_1^{\mathrm{Taylor}}=
    \begin{cases}
    \displaystyle
    \arg\min\limits_{\kappa_1\in\{\kappa_1^{(1)},\kappa_1^{(2)}\}}
    \hat J(\kappa_1), & \Delta>0,\\[1mm]
    \displaystyle {-{B}/{(2A)}}, & \Delta\le0.
    \end{cases}
}
\label{eq:k1_taylor}
\end{equation}
Because the current results are based on the approximations, \(\kappa_1^{\mathrm{Taylor}}\) is further refined by solving the first nonlinear equation in \eqref{eq:exact_system}, i.e., 
\begin{equation}
{
    r_0-\int_0^1
    \cos\left(\sigma(t;\kappa_1,\Gamma-2\kappa_1)\right)dt=0.
}
\label{eq:1d_search_exact}
\end{equation}

\subsection{Solving the Guidance Parameters \(\kappa_1\) and \(\kappa_2\)}
The refined parameter vector
$
    \boldsymbol{\kappa}^{0}
    =
    [\kappa_1^{0},\Gamma-2\kappa_1^{0}]^T
$
is then passed to the nonlinear solver as the initial guess for solving the exact nonlinear system of equations in \eqref{eq:exact_system}. The procedure is summarized in Algorithm~\ref{alg:atac_solver}. 
\begin{algorithm}[htbp]
\caption{Procedure to Solve the Parameters \(\kappa_1\) and \(\kappa_2\)}
\label{alg:atac_solver}
\begin{algorithmic}[1]
\Require Initial condition in \eqref{eq:initial} and desired arrival angle \(\gamma_f\).
\Ensure Exact guidance parameters \(\kappa_1,\kappa_2\).
\Statex \textbf{\textit{Stage 1: Analytical Warm Start}}
\State Compute \(\Gamma\) in \eqref{eq:kappa2_linear}.
\State Evaluate \(A,B,C\) from \eqref{eq:coeff_C} and compute
\(\Delta=B^2-4AC\).
\If{\(\Delta>0\)}
    \State Compute \(\kappa_1^{(1)}\) and \(\kappa_1^{(2)}\).
    \State Evaluate \(\hat J(\kappa_1)\) for both candidates.
    \State Select the candidate with the smaller \(\hat J\).
\Else
    \State Set \(\kappa_1^{\mathrm{Taylor}}\gets -B/(2A)\).
\EndIf
\Statex \textbf{\textit{Stage 2: Solving the Nonlinear System of Equations}}    
\State Obtain \(\kappa_1^{0}\) by solving \eqref{eq:1d_search_exact}.
\State Recover \(\kappa_2^{0}\) using \eqref{eq:kappa2_linear}.
\State Solve \eqref{eq:exact_system} using \([\kappa_1^{0},\kappa_2^{0}]^T\) as the initial guess.
\State \Return \(\kappa_1,\kappa_2\).
\end{algorithmic}
\end{algorithm}

\section{Numerical Simulations}\label{sec:simulations}
This section presents numerical simulations to validate the proposed guidance law. 
MATLAB \textit{fsolve} is used to solve \eqref{eq:1d_search_exact} and \eqref{eq:exact_system}. Unless otherwise stated, the solver tolerances are set to \(10^{-6}\). The physical initial conditions and desired arrival time are normalized such that the arrival time is set to \(t_f=1\), and the resulting normalized solution is scaled back to the physical units. It is worth mentioning that in practice, the UAV is subject to various disturbances and uncertainties, such as wind gusts and sensor noise. Therefore, maintaining a constant-speed flight is challenging. Thanks to our previous work in \cite{wang2026FOV}, which shows that once the guidance law can be generated in real time, the guidance law can be implemented in a closed-loop manner with a sufficiently high update rate to 
deal with different disturbances and uncertainties. Therefore, in this paper, the robustness of the proposed method against disturbances and uncertainties is not explicitly analyzed. Instead, the focus is on demonstrating the computational efficiency and solution quality of the proposed method in closed-loop implementation and comparative studies.
\subsection{Closed-Loop Implementation}
The initial conditions are set to \(r_0=5000~\mathrm{m}\), \(\lambda_0=0^\circ\), and \(\gamma_0=30^\circ\). The desired arrival angle is \(\gamma_f=-60^\circ\), the desired arrival time is \(t_f=35~\mathrm{s}\), and the UAV speed is \(V=200~\mathrm{m/s}\). The guidance update interval is \(\Delta t=0.01~\mathrm{s}\). The parameters \(\kappa_1\) and \(\kappa_2\) are recomputed based on state feedback at each guidance update. The analytical initialization is required only in the first guidance cycle; in subsequent cycles, the converged parameters from the previous cycle are used as the initial guess. The acceleration command is saturated by \(\pm3g\), and the guidance law switches to PN guidance with a gain of \(3\) when the relative range drops below 100 m. The simulation terminates once the range reaches \(1~\mathrm{m}\).

Figure~\ref{fig:sim_results} shows the closed-loop results. The UAV reaches the destination at the prescribed time and with the desired arrival angle while satisfying the acceleration constraint. The guidance parameters vary smoothly and approach zero in the terminal phase, which is consistent with the fact that the range, look angle, and remaining flight time all vanish at arrival. The nonlinear solver requires the largest number of iterations of ten only in the first guidance cycle; afterward, the previous-cycle solution provides an excellent initial guess, and the solver typically converges within one iteration or requires no additional iteration. This confirms the real-time computational efficiency of the proposed method.
\begin{figure}[t]
    \centering

    \subfloat[Flight trajectory]{
        \includegraphics[width=0.47\textwidth]{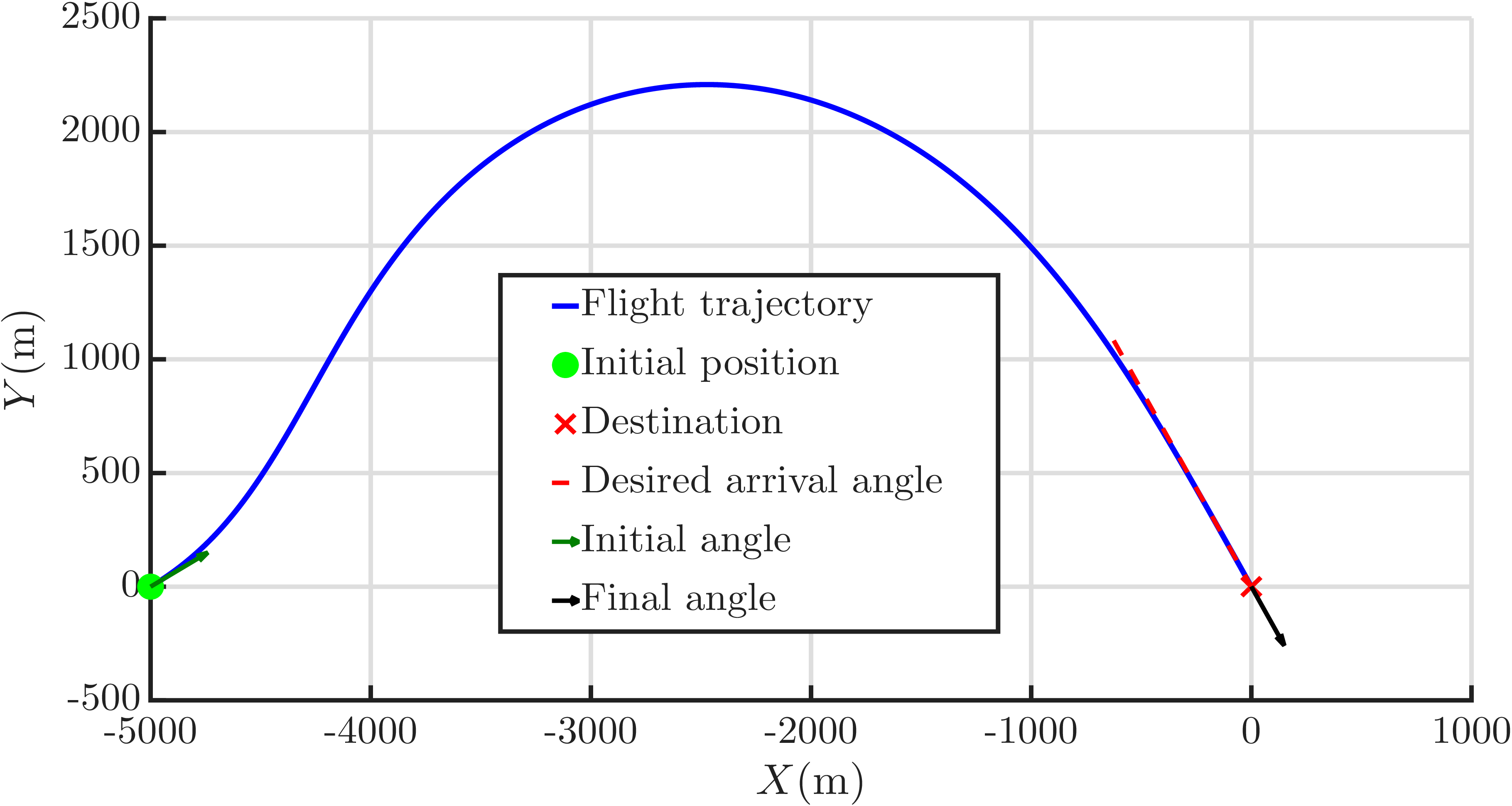}
        \label{fig:closed_loop_traj}
    }
    \hfill
    \subfloat[Guidance parameters \(\kappa_1\) and \(\kappa_2\)]{
        \includegraphics[width=0.47\textwidth]{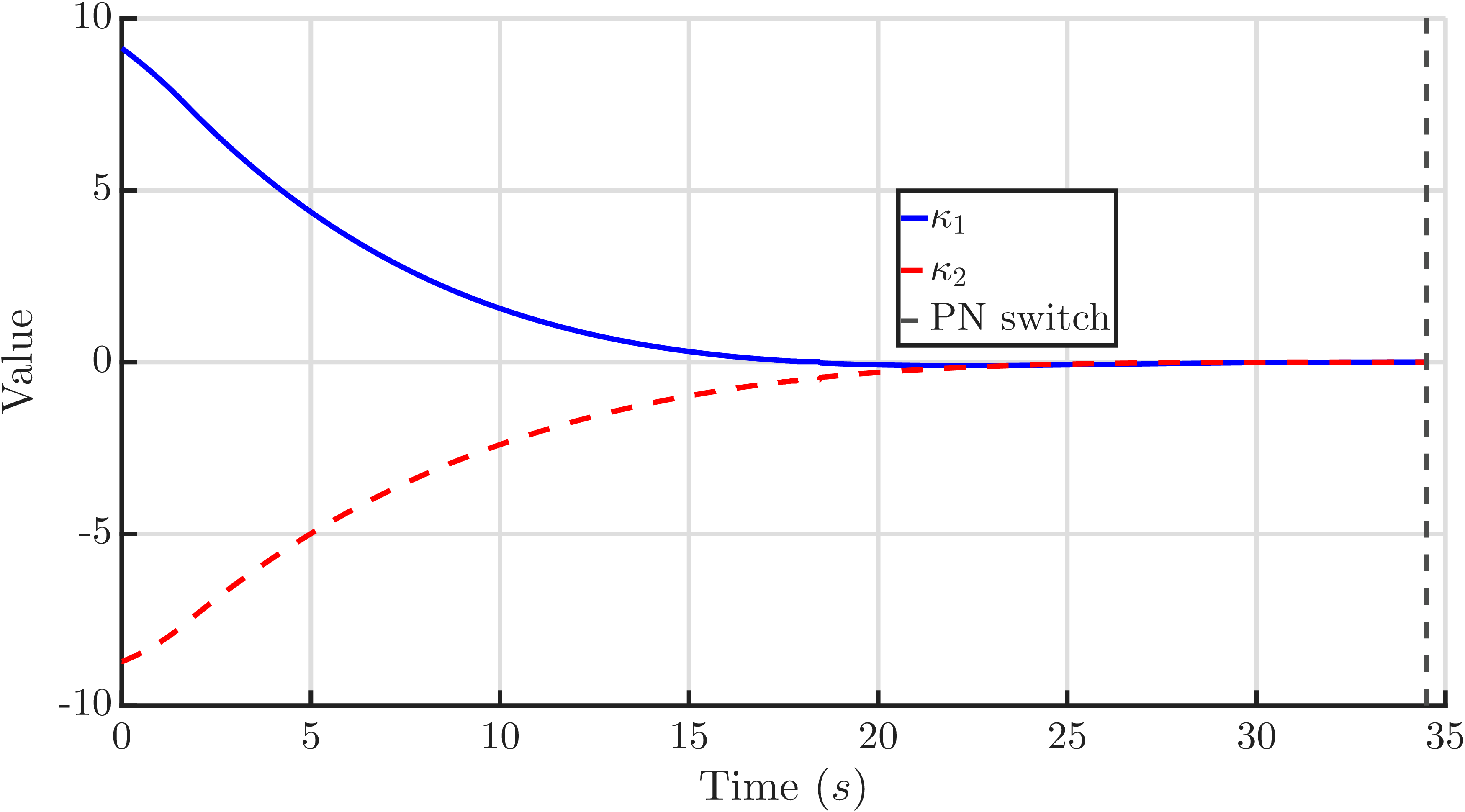}
        \label{fig:fig_guidanceP}
    }

    \vspace{0.2cm}

    \subfloat[Acceleration and angles]{
        \includegraphics[width=0.47\textwidth]{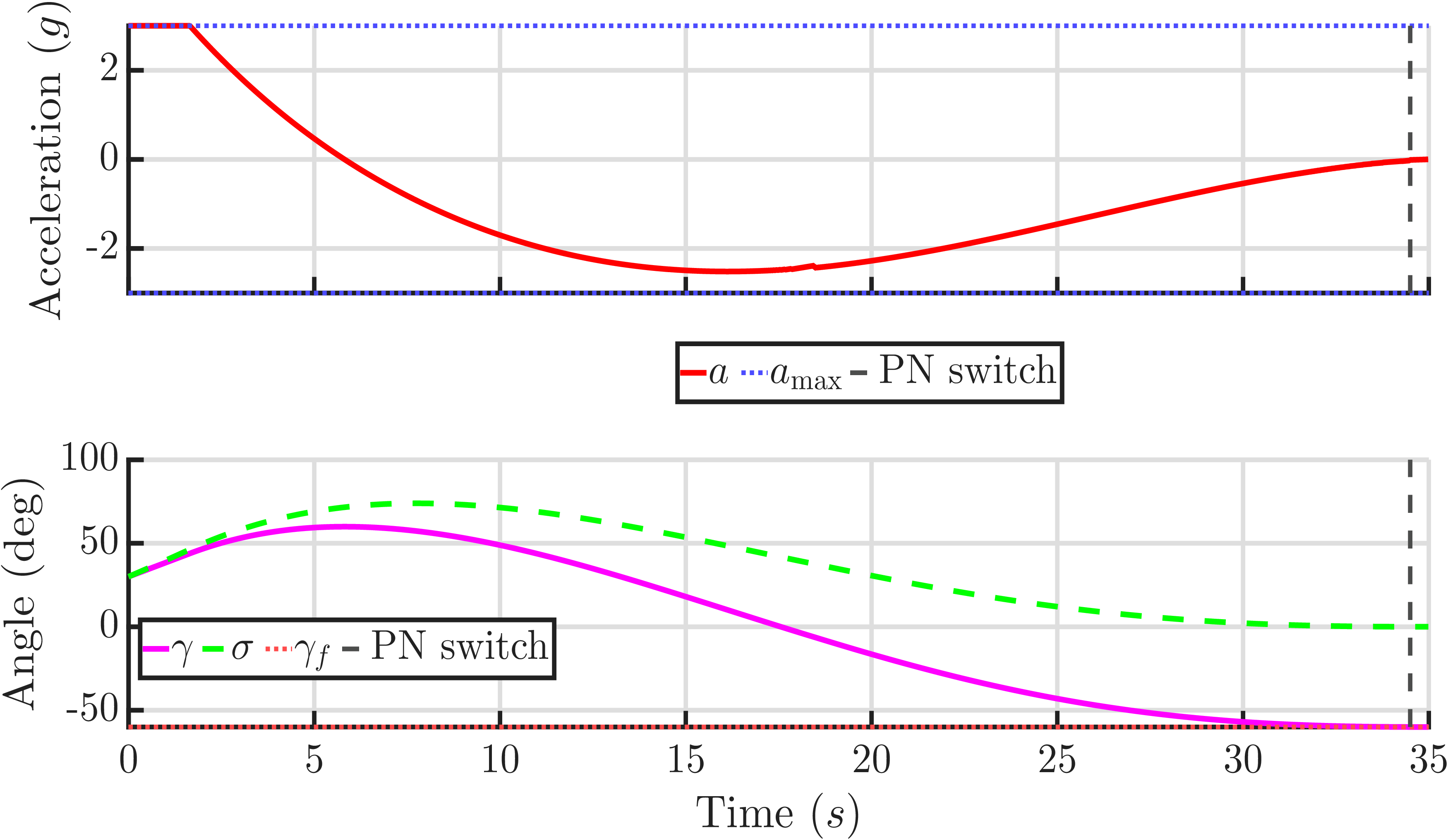}
        \label{fig:angles}
    }
    \hfill
    \subfloat[Number of iterations]{
        \includegraphics[width=0.47\textwidth]{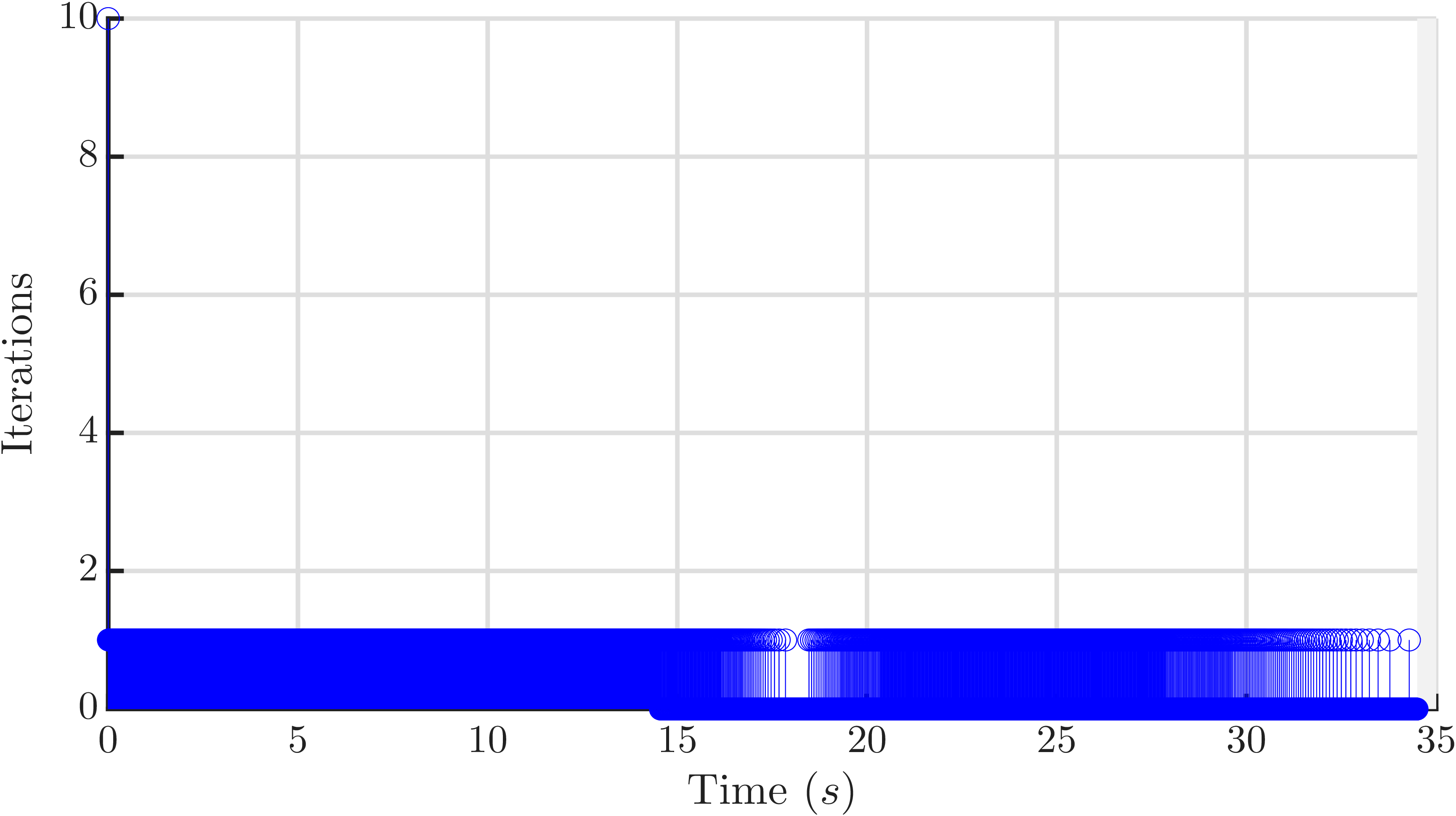}
        \label{fig:iterations}
    }

    \caption{Simulation results of the proposed method in closed-loop.}
    \label{fig:sim_results}
\end{figure}

\subsection{Comparative Study}
The proposed method is compared with three alternatives. The first is the trajectory-shaping method in \cite{kang2019generalized}. The second is a geometry-based virtual-target guidance law proposed by Hou et al. in \cite{hou2023optimal}. The third is an open-loop optimal control solver GPOPS-II from \cite{patterson2014gpops}, which serves as a benchmark. Although GPOPS-II can generate high-quality solutions, it usually requires a good initial guess and is computationally expensive, making it not directly suitable for onboard guidance. Two cases will be considered to demonstrate the performance of the proposed method and the alternatives. 
\subsubsection{Case A}
We consider an initial condition from \cite{hou2023optimal}, i.e.,
\(r_0=15000~\mathrm{m}\), \(\lambda_0=90^\circ\), \(\gamma_0=0^\circ\), \(\gamma_f=-60^\circ\), and \(t_f=90~\mathrm{s}\). The UAV speed is \(V=250~\mathrm{m/s}\). The results are shown in Fig.~\ref{fig:comparative_results}. The method in \cite{kang2019generalized} produces the same trajectory as the proposed method in this case because both approaches employ an equivalent fourth-order polynomial and enforce the same boundary conditions; therefore, only the proposed method is shown for clarity. All methods reach the destination with the desired arrival angle at the desired arrival time. However, the results from the proposed method are closer to those from GPOPS-II, as observed from the flight trajectory and guidance command profiles from Figs.~\ref{fig:comparative_trajectory}--\ref{fig:comparative_command}. The control effort is evaluated as \(J=\int_0^{t_f}a^2(t)\,dt\), and the results are given in Table~\ref{table_control_effort}. It is clear that the proposed method achieves a control effort nearly identical to the GPOPS-II benchmark, while the method in \cite{hou2023optimal} yields a slightly higher control effort.
\begin{figure}[t]
    \centering

    \subfloat[Flight trajectories]{
        \includegraphics[width=0.4\textwidth]{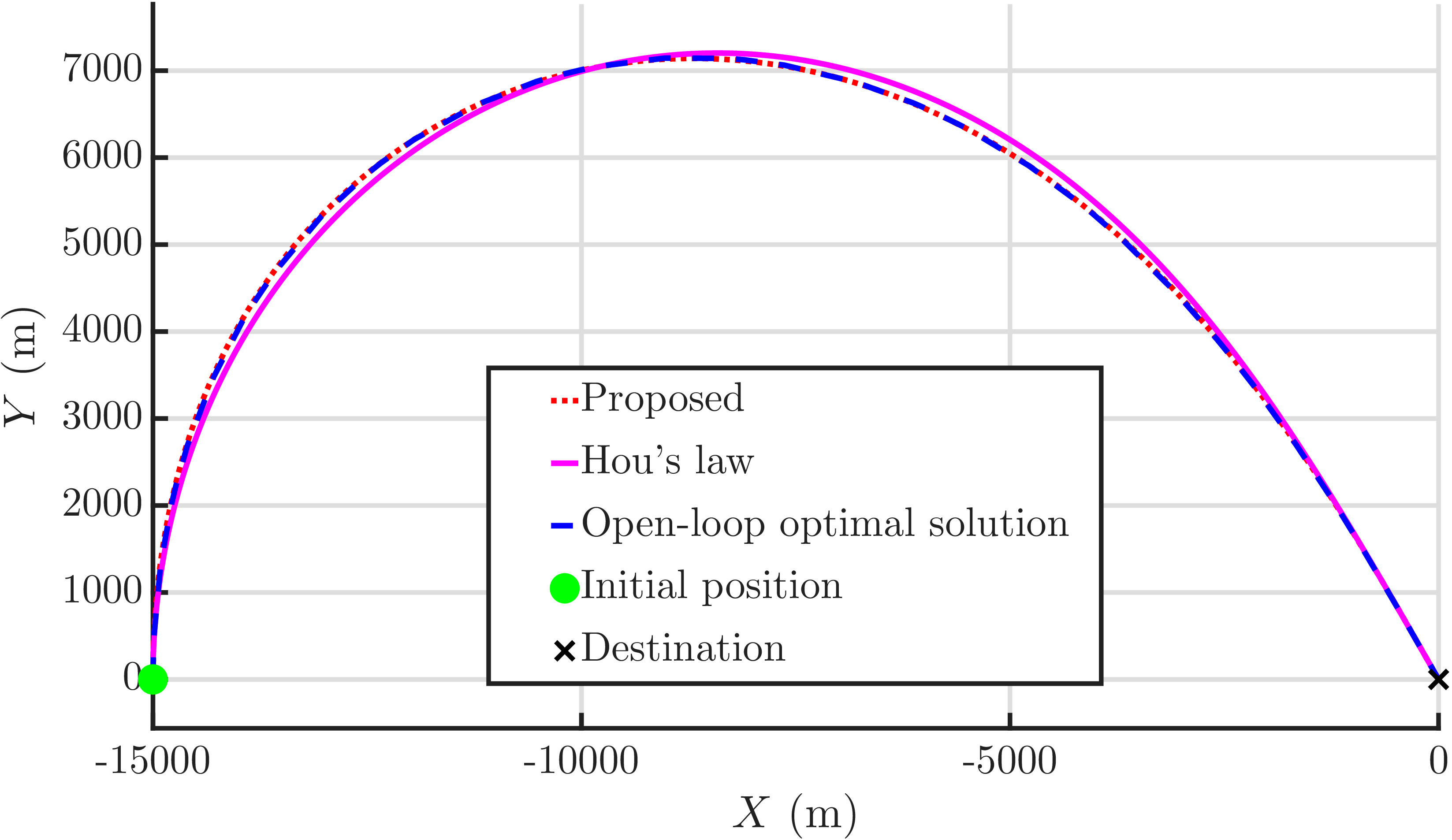}
        \label{fig:comparative_trajectory}
    }

    \vspace{0.2cm}

    \subfloat[Flight path angles]{
        \includegraphics[width=0.4\textwidth]{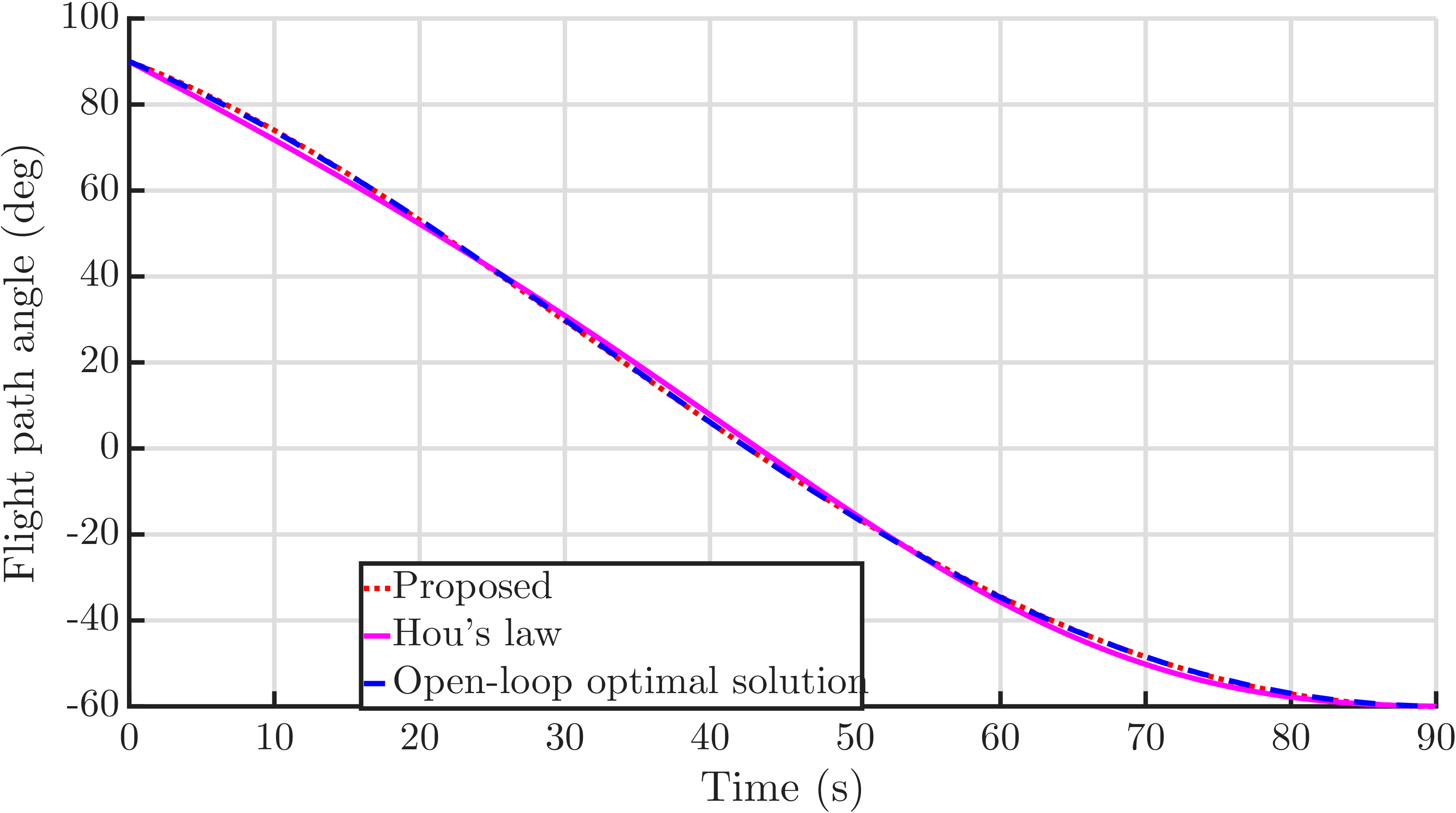}
        \label{fig:comparative_angle}
    }

    \vspace{0.2cm}

    \subfloat[Guidance commands]{
        \includegraphics[width=0.4\textwidth]{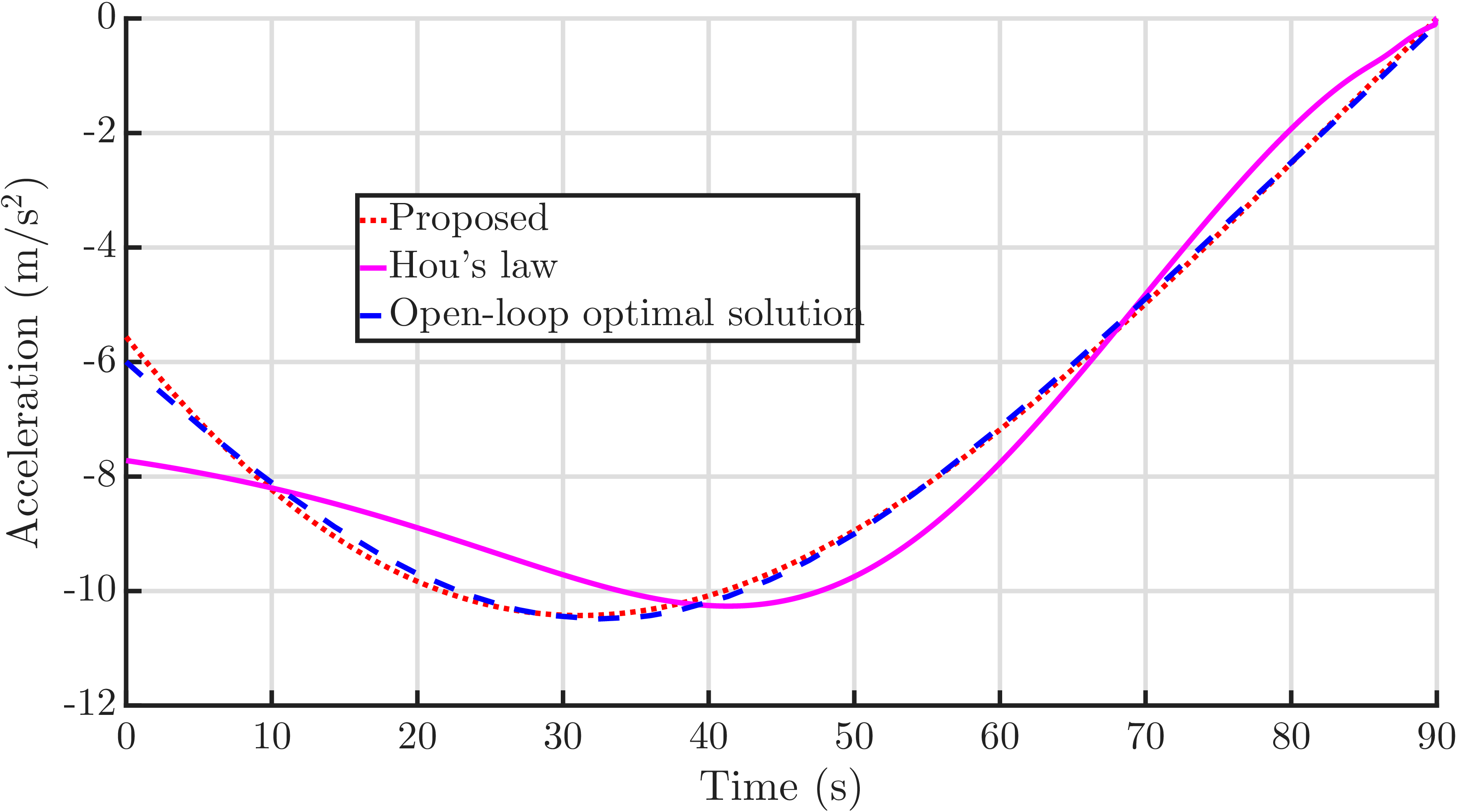}
        \label{fig:comparative_command}
    }

    \caption{Comparative study for Case A.}
    \label{fig:comparative_results}
\end{figure}
\begin{table}[t]
\caption{Comparison of control effort for Case A}
\label{table_control_effort}
\centering
\begin{tabular}{lc}
\toprule
Method & \(J=\int_{0}^{t_f}a^2(t)\,dt~(\mathrm{m^2\,s^{-3}})\) \\
\midrule
Proposed  & \(5.558\times10^{3}\) \\
Hou's law & \(5.589\times10^{3}\) \\
GPOPS-II  & \(5.555\times10^{3}\) \\
\bottomrule
\end{tabular}
\end{table}
\subsubsection{Case B}
We next consider a highly nonlinear case with \(r_0=10000~\mathrm{m}\), \(\lambda_0=0^\circ\), \(\gamma_0=120^\circ\), \(\gamma_f=120^\circ\), and \(t_f=150~\mathrm{s}\). The UAV speed is \(V=200~\mathrm{m/s}\). Both the initial and final flight-path angles are larger than \(90^\circ\), and the desired arrival time is long, which leads to a turn-back-type trajectory and makes the problem particularly challenging. Among the considered methods, the proposed method and GPOPS-II successfully find solutions that satisfy all the constraints, as shown in Fig.~\ref{fig:comparative_results_2}. The method in \cite{kang2019generalized} fails to obtain a feasible solution for this case, which highlights the importance of the proposed initialization strategy. Hou's law is efficient in less challenging cases, but it becomes less reliable when the engagement scenario is highly nonlinear.

From Fig.~\ref{fig:comparative_trajectory_2}, the proposed method yields a trajectory very close to that from GPOPS-II. In Fig.~\ref{fig:comparative_command_2}, the analytical approximation \(\hat a(t)\) in \eqref{eq:accel_linear} closely matches the actual acceleration \(a(t)\), validating the usefulness of the approximations in \eqref{eq:approximations}. Moreover, the actual acceleration command of the proposed method is also close to the open-loop optimal solution from
GPOPS-II. In terms of control effort, the proposed method yields \(J=1.7876\times10^4~\mathrm{m^2\,s^{-3}}\), which is close to the GPOPS-II benchmark value \(J=1.7634\times10^4~\mathrm{m^2\,s^{-3}}\).
\begin{figure}[t]
    \centering
    \subfloat[Flight trajectories from the proposed method and GPOPS-II]{
        \includegraphics[width=0.47\textwidth]{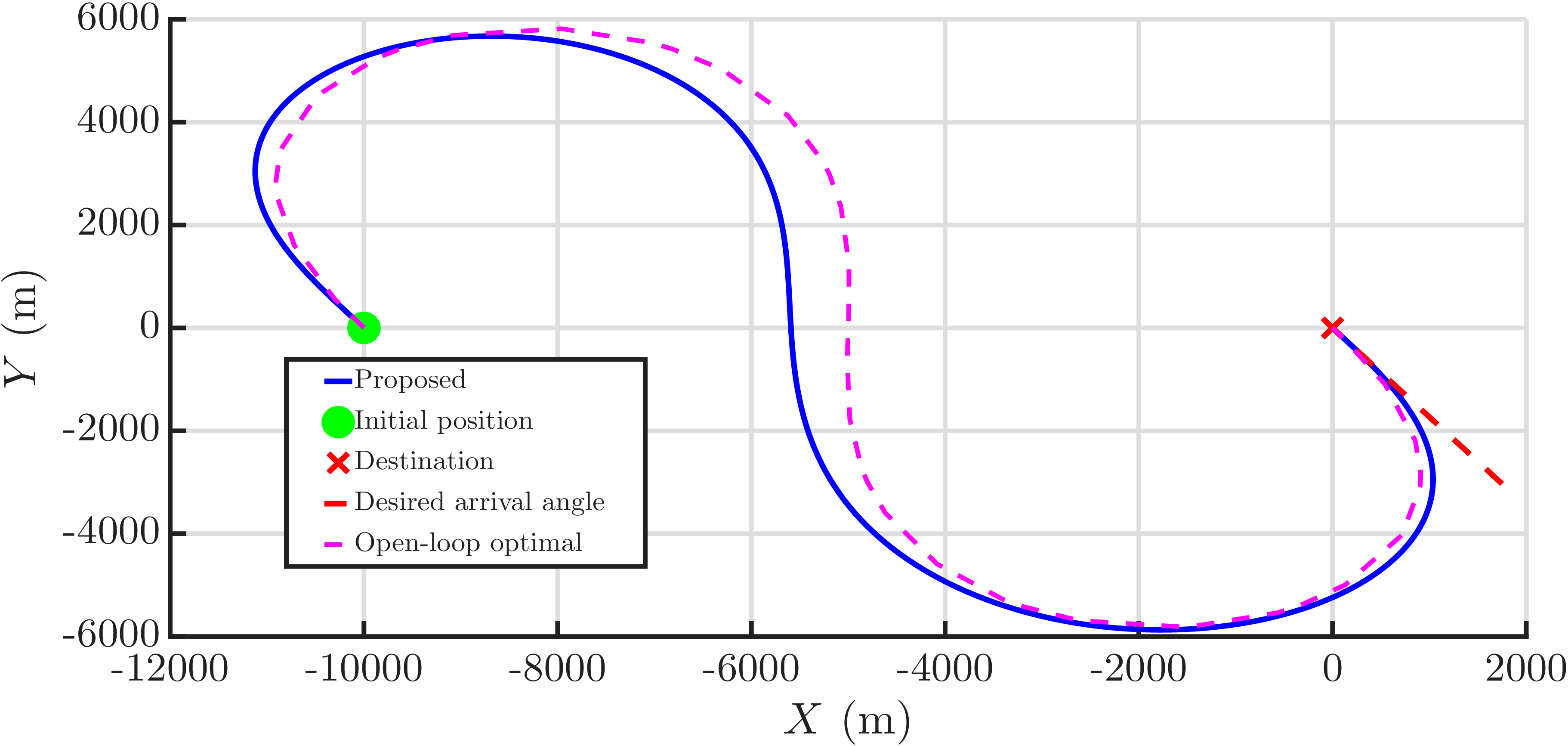}
        \label{fig:comparative_trajectory_2}
    }
    \hfill
    \subfloat[Angles from the proposed method]{
        \includegraphics[width=0.47\textwidth]{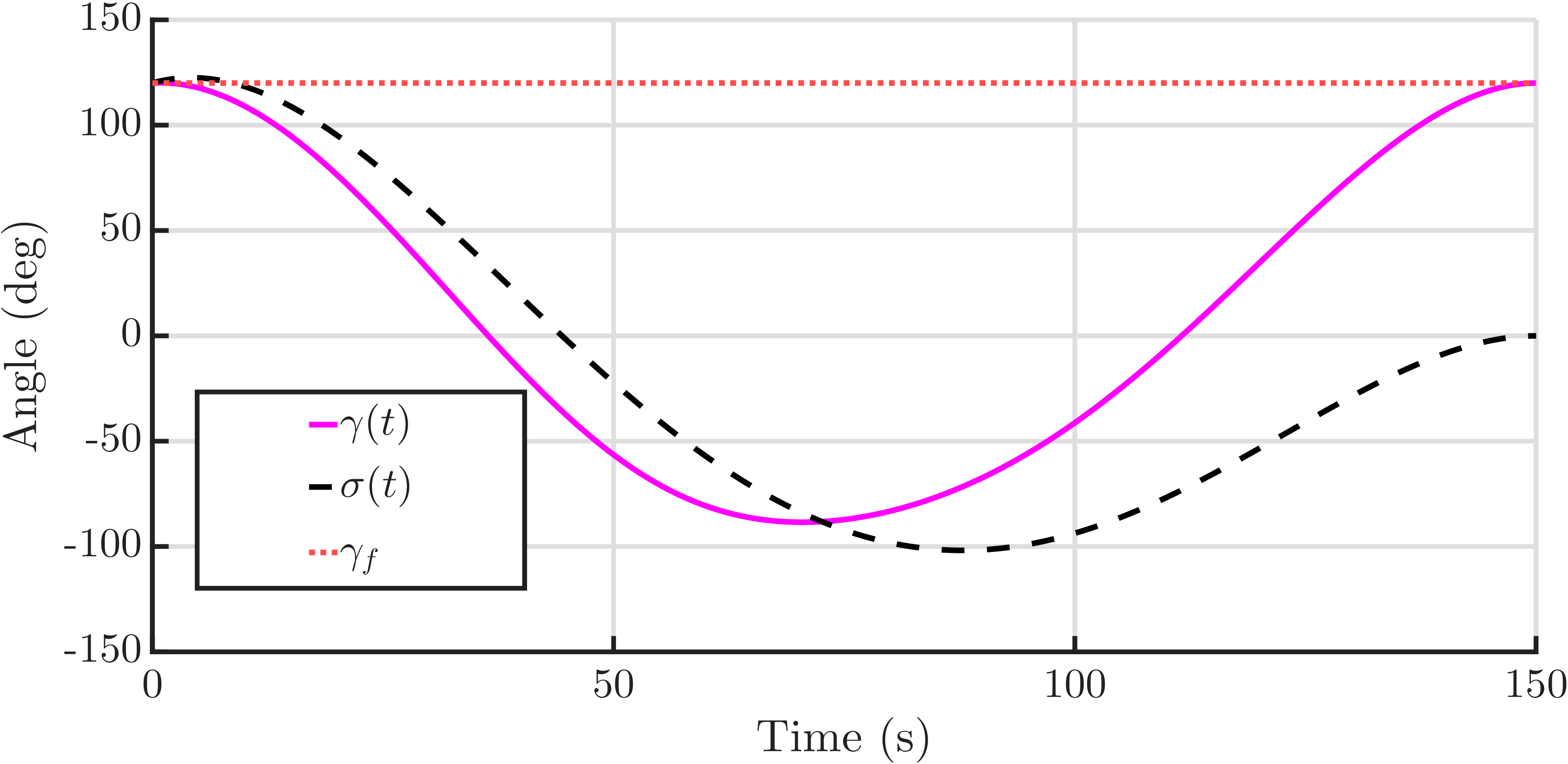}
        \label{fig:comparative_angle_2}
    }

    \vspace{0.2cm}

    \subfloat[Guidance from the proposed method and GPOPS-II]{
        \includegraphics[width=0.47\textwidth]{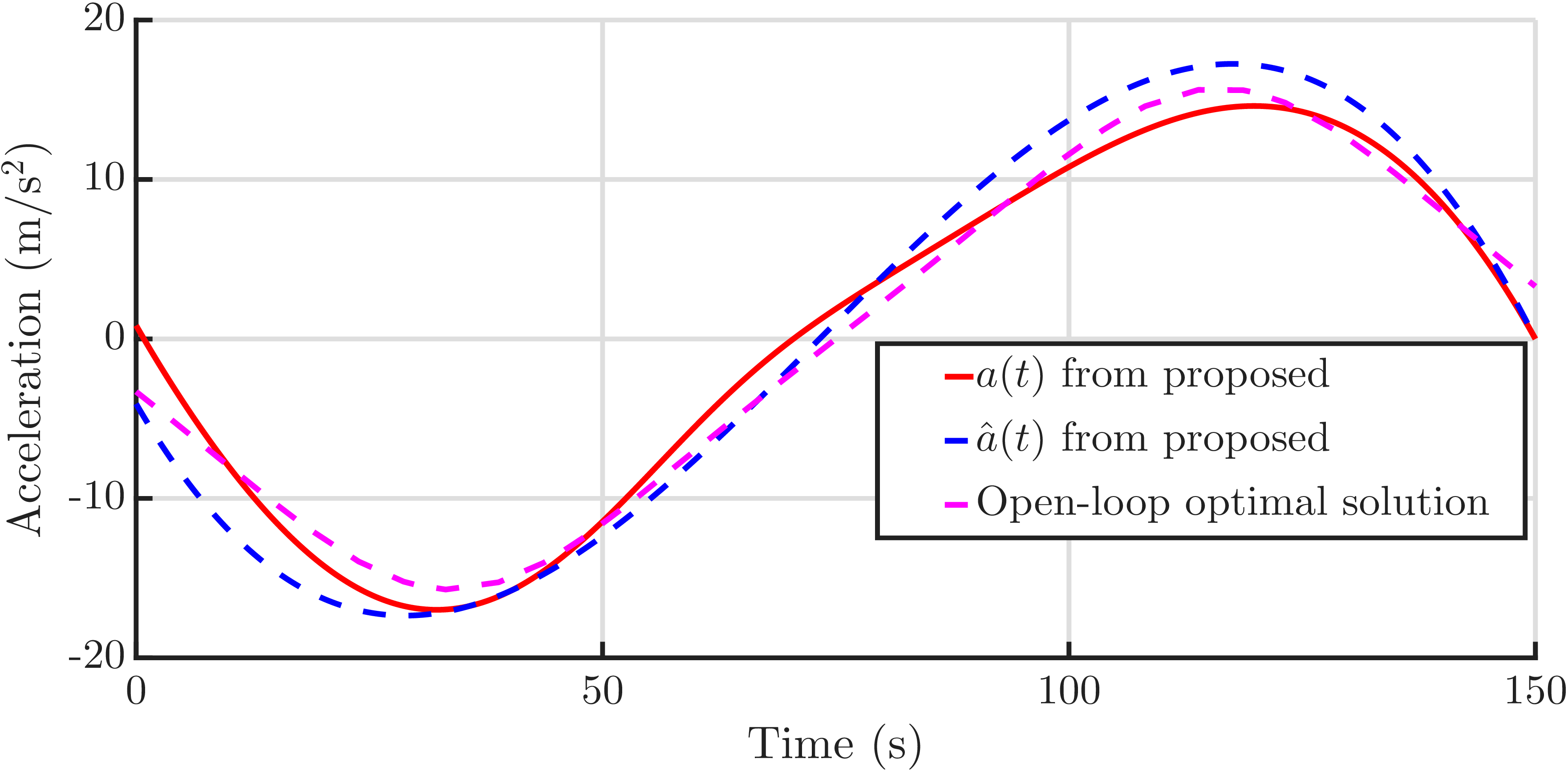}
        \label{fig:comparative_command_2}
    }
    \hfill
    \subfloat[Range from the proposed method and GPOPS-II]{
        \includegraphics[width=0.47\textwidth]{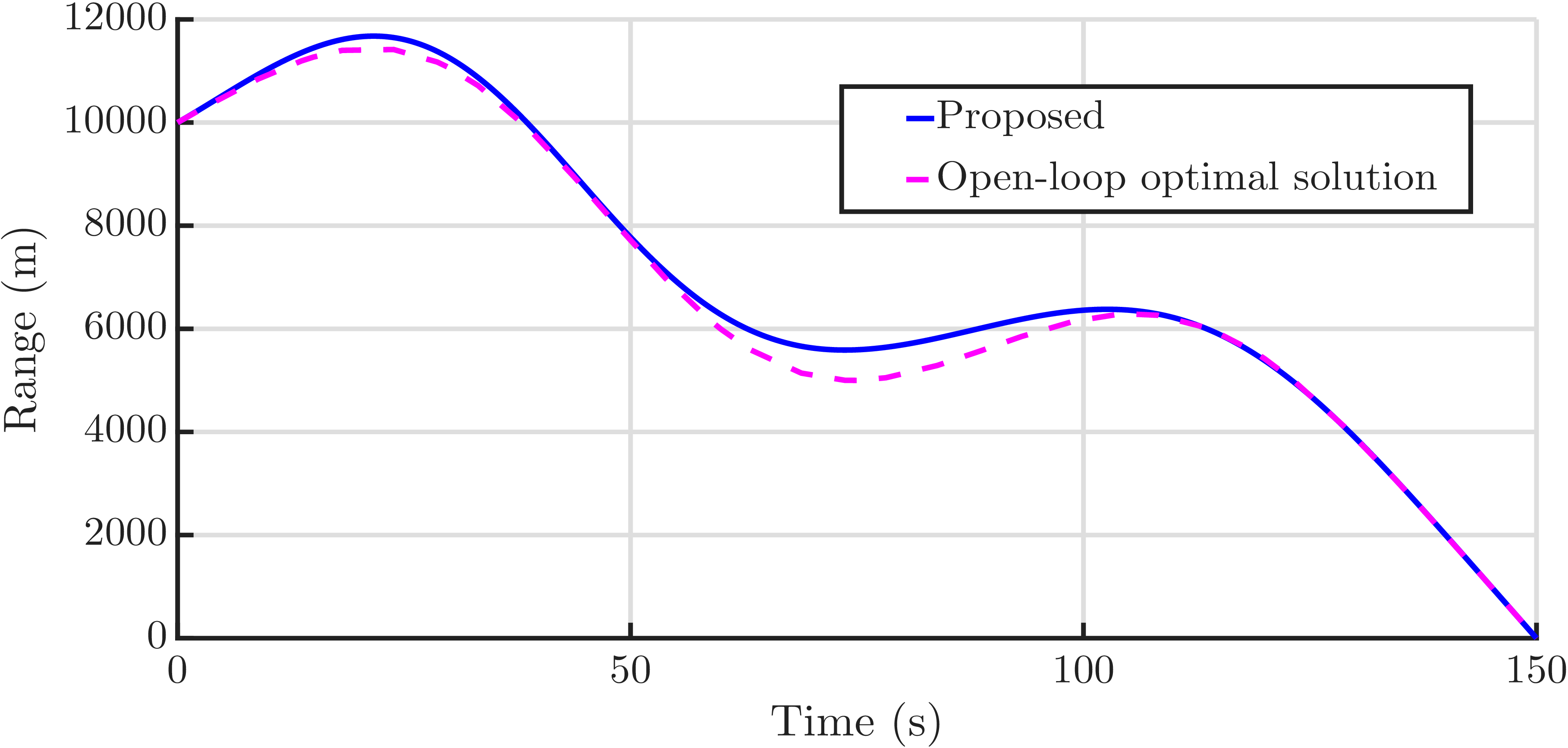}
        \label{fig:comparative_range_2}
    }

    \caption{Comparative study for Case B.}
    \label{fig:comparative_results_2}
\end{figure}

\subsection{Application to Formation Establishment}
We further apply the proposed guidance law to a formation-establishment scenario, where three UAVs are required to arrive at their respective destinations simultaneously and form a prescribed triangular pattern. The desired destinations are located at the vertices of an equilateral triangle with a side length of \(1000~\mathrm{m}\), centered at the origin. The desired terminal flight-path angles of all three UAVs are set to \(\gamma_f=0^\circ\), so that the vehicles share the same final heading direction after the formation is established. The common desired arrival time is selected as \(t_f=55~\mathrm{s}\).
The initial positions and flight path angles of the three UAVs are given by
\begin{align*}
    \mathrm{UAV~1:}\quad &(x_0,y_0,\gamma_0)
    =(-8.5~\mathrm{km},\,0.35~\mathrm{km},\,2^\circ),\\
    \mathrm{UAV~2:}\quad &(x_0,y_0,\gamma_0)
    =(-9.5~\mathrm{km},\,-0.70~\mathrm{km},\,2.5^\circ),\\
    \mathrm{UAV~3:}\quad &(x_0,y_0,\gamma_0)
    =(-10.5~\mathrm{km},\,-1.10~\mathrm{km},\,4^\circ).
\end{align*}
The flight speeds of UAV~1, UAV~2, and UAV~3 are set to \(160~\mathrm{m/s}\), \(175~\mathrm{m/s}\), and \(190~\mathrm{m/s}\), respectively. 

The simulation results are shown in Fig.~\ref{fig:cooperative_results}. The trajectory profiles show that the three UAVs reach the three prescribed vertices of the equilateral triangle at the same arrival time. The terminal direction arrows indicate that all vehicles arrive with the common desired flight path angle, thereby satisfying the formation-orientation requirement.  In addition, the flight path angle profiles converge smoothly to \(\gamma_f=0^\circ\). These results verify that the proposed guidance law can coordinate multiple UAVs in a formation-establishment task. Since the terminal accelerations are zero, the established formation can be maintained without any additional maneuver.
\begin{figure}[t]
    \centering

    \subfloat[Flight trajectories]{
        \includegraphics[width=0.4\textwidth]{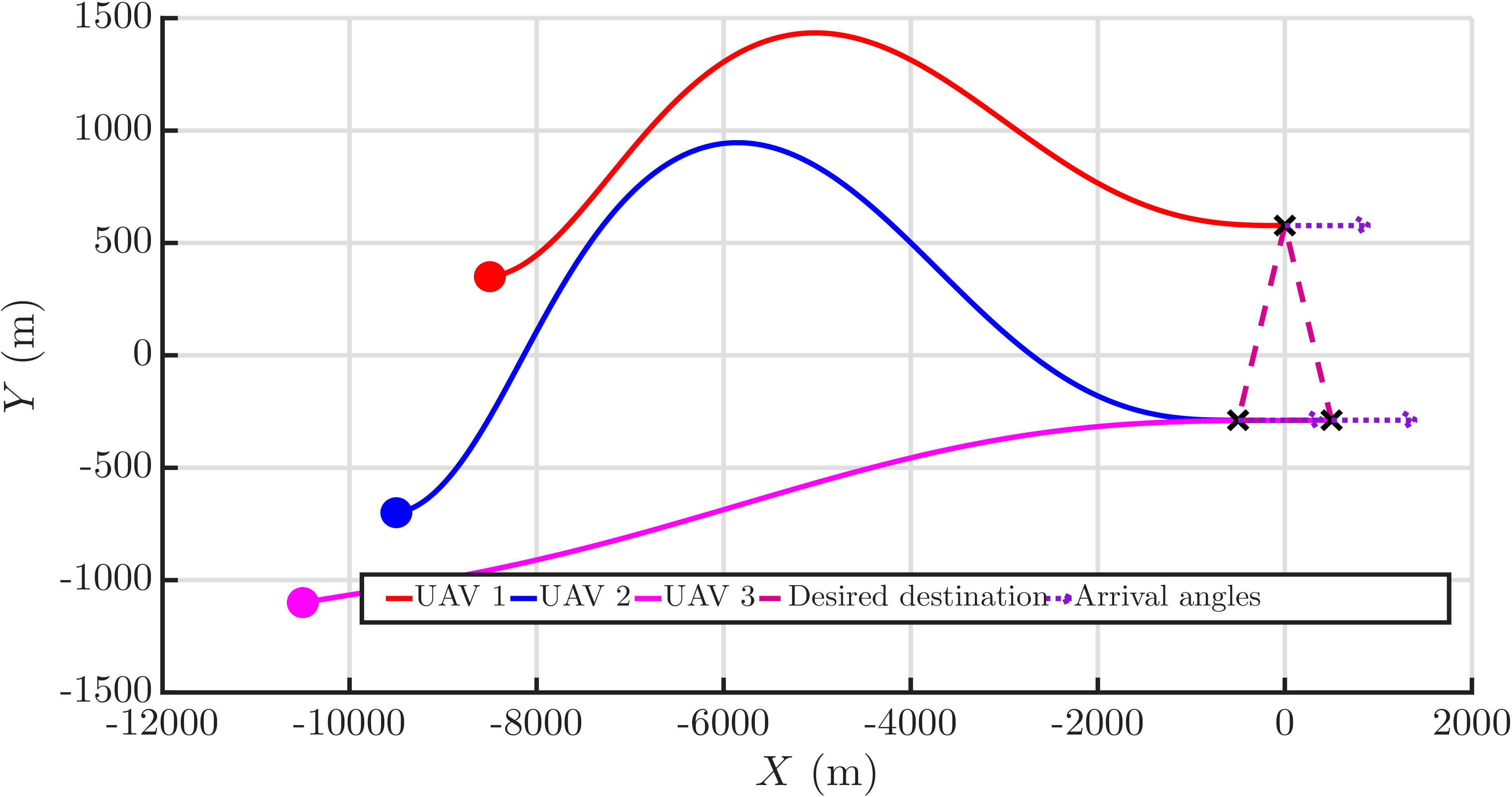}
        \label{fig:formation_traj}
    }

    \vspace{0.2cm}

    \subfloat[Guidance command profiles]{
        \includegraphics[width=0.4\textwidth]{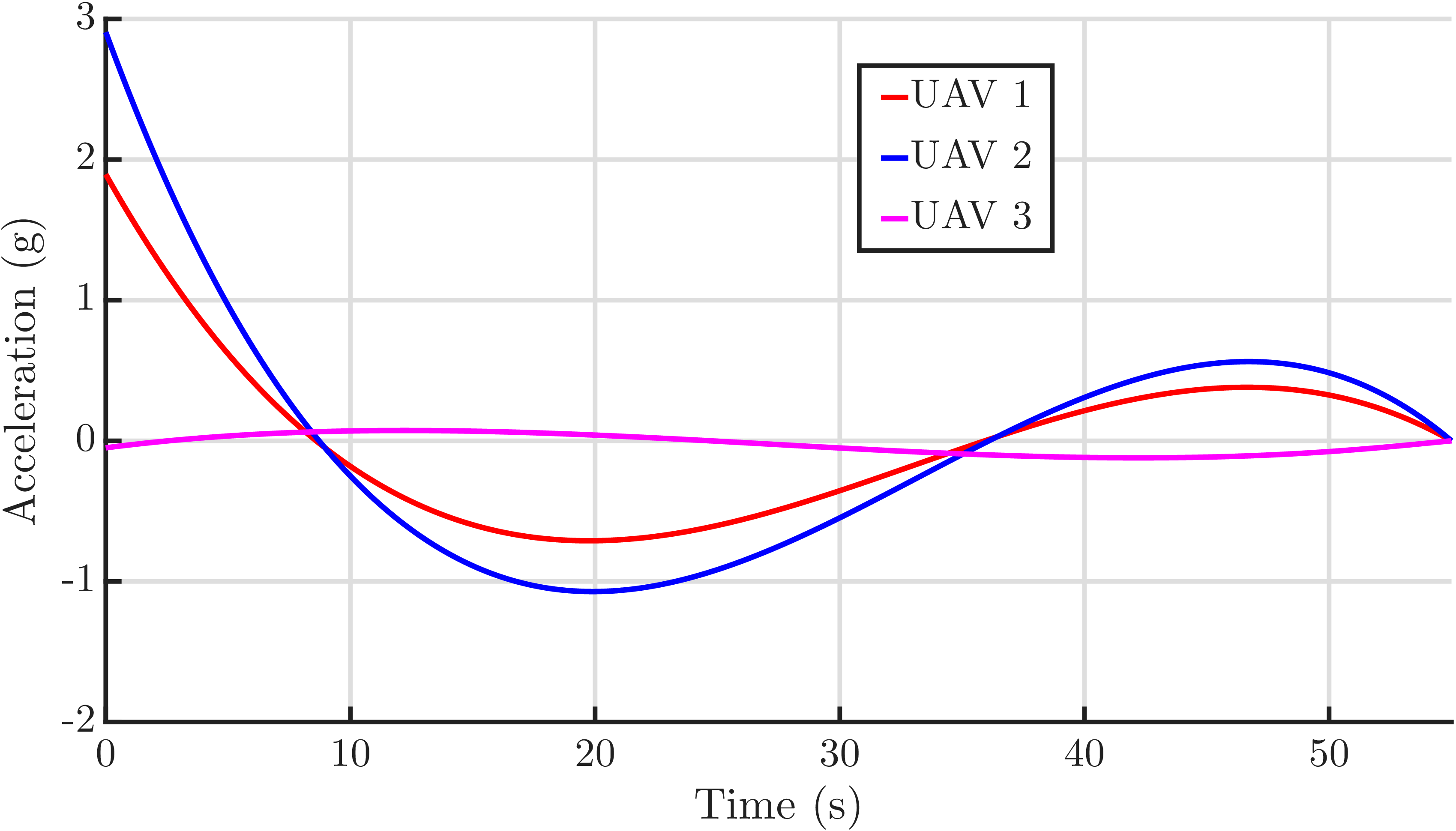}
        \label{fig:formation_acc}
    }

    \vspace{0.2cm}

    \subfloat[Flight path angle profiles]{
        \includegraphics[width=0.4\textwidth]{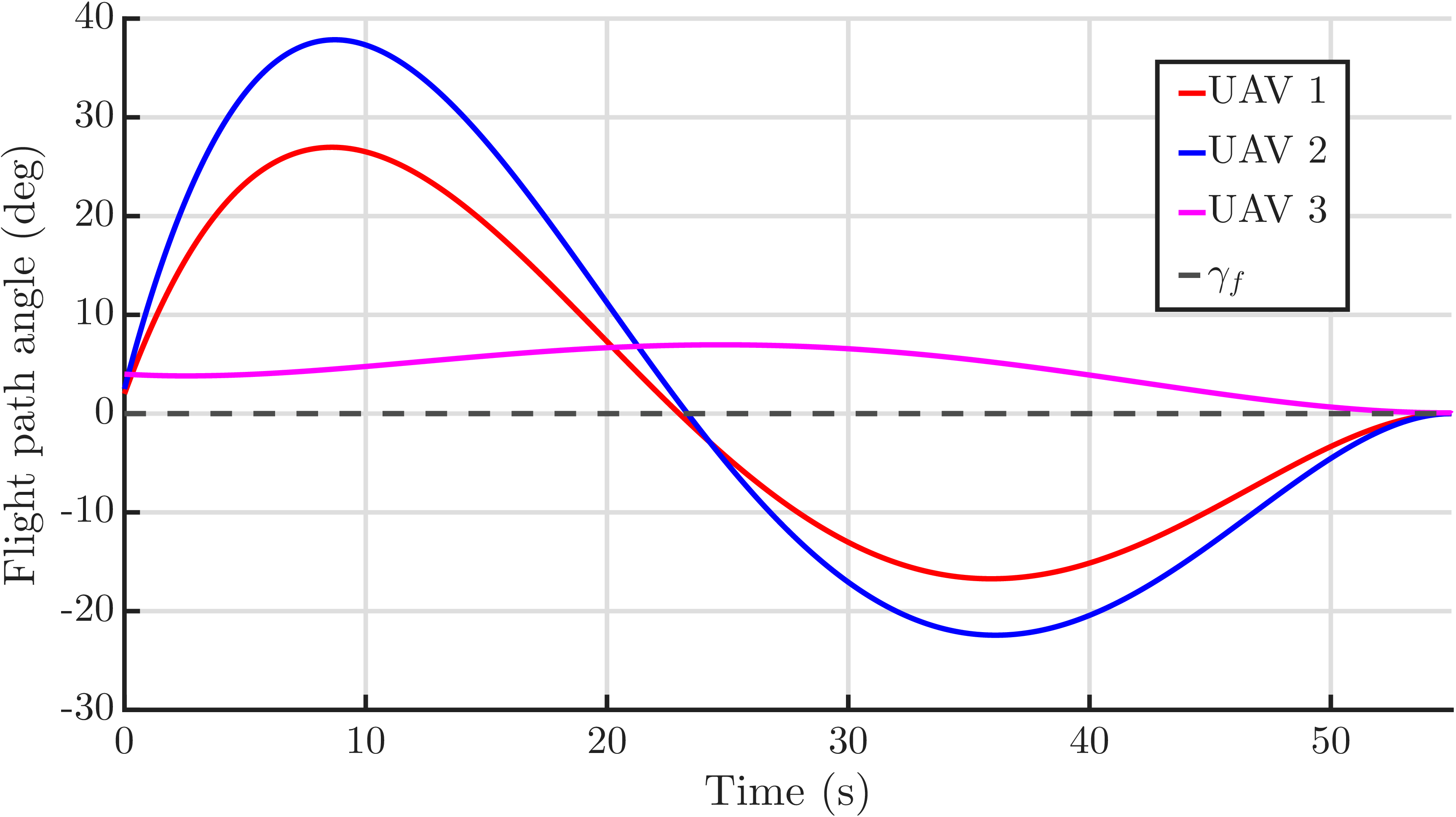}
        \label{fig:formation_angle}
    }

    \caption{Simulation results of the formation establishment for three UAVs.}
    \label{fig:cooperative_results}
\end{figure}
\section{Conclusions}\label{sec:conclusions}
This paper presented a nonlinear trajectory shaping guidance law for arrival-time-and-arrival-angle control of a constant-speed fixed-wing UAV. By parameterizing the look angle with a fourth-order polynomial, the terminal look angle and look angle rate constraints are naturally satisfied. In this way, the highly nonlinear guidance problem is reduced to solving for two shaping parameters.
A two-stage solution procedure is developed to compute the shaping parameters efficiently. The analytical warm start provides a reliable initial guess, and the subsequent nonlinear solve enforces the exact terminal conditions. Numerical simulations demonstrated that the proposed method achieves fast convergence in closed-loop implementation and produces feasible solutions very close to the open-loop optimal solution benchmark. In addition, the proposed method is applicable to highly nonlinear scenarios where some existing methods fail to find feasible solutions in the test case.

\bibliography{References}

@article{wang2026FOV,
  title={Look-Angle-Shaped Impact Time Control Guidance
with Field-of-View Constraints},
  author={Wang, Kun and Wei, Zhenyu and Wang, Pengyu},
  journal={Aerospace Science and Technology},
  year={2026},
  publisher={Elsevier},
  pages = {112131},
  doi = {10.1016/j.ast.2026.112131}
}

@article{singh2025terminal,
  title={Terminal Time Constrained Guidance Against Stationary Target with Bounded Input and Rate},
  author={Singh, Swati and Kumar, Shashi Ranjan and Mukherjee, Dwaipayan},
  journal={IEEE Control Systems Letters},
  year={2025},
  publisher={IEEE},
  doi = {10.1109/LCSYS.2025.3578569}
}

@article{kumar2018sliding,
  title={Sliding mode guidance for impact time and angle constraints},
  author={Kumar, Shashi Ranjan and Ghose, Debasish},
  journal={Proceedings of the Institution of Mechanical Engineers, Part G: Journal of Aerospace Engineering},
  volume={232},
  number={16},
  pages={2961--2977},
  year={2018},
  publisher={SAGE Publications Sage UK: London, England},
  doi = {10.1177/0954410017719304}
}

@article{kang2019generalized,
  title={Generalized impact time and angle control via look-angle shaping},
  author={Kang, Shen and Tekin, Raziye and Holzapfel, Florian},
  journal={Journal of Guidance, Control, and Dynamics},
  volume={42},
  number={3},
  pages={695--702},
  year={2019},
  publisher={American Institute of Aeronautics and Astronautics},
  doi = {10.2514/1.G003765}
}

@inproceedings{wang2026trajectory,
  title     = {Trajectory Shaping Guidance for Field-of-View Constrained Impact Angle Control},
  author    = {Wang, K.},
  booktitle = {The CEAS Conference on Guidance, Navigation and Control, EuroGNC 2026},
  publisher = {Council of European Aerospace Societies (CEAS)},
  year      = {2026},
  doi       = {10.82124/CEAS-GNC-2026-027}
}

@article{majumder2023three,
  title={Three-dimensional nonlinear impact time guidance considering field-of-view constraints},
  author={Majumder, Kakoli and Kumar, Shashi Ranjan},
  journal={IEEE Control Systems Letters},
  volume={7},
  pages={3848--3853},
  year={2023},
  publisher={IEEE},
  doi = {10.1109/LCSYS.2023.3341447}
}

@article{patterson2014gpops,
  title={{GPOPS-II}: A {MATLAB} software for solving multiple-phase optimal control problems using hp-adaptive Gaussian quadrature collocation methods and sparse nonlinear programming},
  author={Patterson, Michael A and Rao, Anil V},
  journal={ACM Transactions on Mathematical Software (TOMS)},
  volume={41},
  number={1},
  pages={1--37},
  year={2014},
  publisher={ACM New York, NY, USA},
  doi = {10.1145/2558904}
}

@article{wang2025nonlinear,
  title={Nonlinear Optimal Impact Angle Control Guidance With Acceleration Constraints},
  author={Wang, Kun and Lu, Fangmin and Chen, Zheng},
  journal={IEEE Transactions on Aerospace and Electronic Systems},
  volume={61},
  number={4},
  pages={8907--8921},
  year={2025},
  doi = {10.1109/TAES.2025.3551283}
}

@article{kim2013augmented,
  title={Augmented polynomial guidance with impact time and angle constraints},
  author={Kim, Tae-Hun and Lee, Chang-Hun and Jeon, In-Soo and Tahk, Min-Jea},
  journal={IEEE Transactions on Aerospace and Electronic Systems},
  volume={49},
  number={4},
  pages={2806--2817},
  year={2013},
  publisher={IEEE},
  doi = {10.1109/TAES.2013.6621856}
}

@article{zhang2013guidance,
  title={Guidance law with impact time and impact angle constraints},
  author={Zhang, Youan and Ma, Guoxin and Liu, Aili},
  journal={Chinese Journal of Aeronautics},
  volume={26},
  number={4},
  pages={960--966},
  year={2013},
  publisher={Elsevier},
  doi = {10.1016/j.cja.2013.04.037}
}

@article{chen2019optimal,
  title={Optimal control based guidance law to control both impact time and impact angle},
  author={Chen, Xiaotian and Wang, Jinzhi},
  journal={Aerospace Science and Technology},
  volume={84},
  pages={454--463},
  year={2019},
  publisher={Elsevier},
  doi = {10.1016/j.ast.2018.10.036}
}

@article{wang2022nonlinear,
  title={Nonlinear optimal guidance for intercepting stationary targets with impact-time constraints},
  author={Wang, Kun and Chen, Zheng and Wang, Han and Li, Jun and Shao, Xueming},
  journal={Journal of Guidance, Control, and Dynamics},
  volume={45},
  number={9},
  pages={1614--1626},
  year={2022},
  publisher={American Institute of Aeronautics and Astronautics},
  doi = {10.2514/1.G006666}
}

@article{chen2023elongation,
  title={Elongation of curvature-bounded path},
  author={Chen, Zheng and Wang, Kun and Shi, Heng},
  journal={Automatica},
  volume={151},
  pages={110936},
  year={2023},
  publisher={Elsevier},
  doi = {10.1016/j.automatica.2023.110936}
}

@article{shanmugavel2010co,
  title={Co-operative path planning of multiple UAVs using Dubins paths with clothoid arcs},
  author={Shanmugavel, Madhavan and Tsourdos, Antonios and White, Brian and {\.Z}bikowski, Rafa{\l}},
  journal={Control Engineering Practice},
  volume={18},
  number={9},
  pages={1084--1092},
  year={2010},
  publisher={Elsevier},
  doi = {10.1016/j.conengprac.2009.02.010}
}

@article{lee2007guidance,
  title={Guidance law to control impact time and angle},
  author={Lee, Jin-Ik and Jeon, In-Soo and Tahk, Min-Jea},
  journal={IEEE Transactions on Aerospace and Electronic Systems},
  volume={43},
  number={1},
  pages={301--310},
  year={2007},
  publisher={IEEE},
  doi = {10.1109/TAES.2007.357135}
}

@article{shaferman2008linear,
  title={Linear quadratic guidance laws for imposing a terminal intercept angle},
  author={Shaferman, Vitaly and Shima, Tal},
  journal={Journal of Guidance, Control, and Dynamics},
  volume={31},
  number={5},
  pages={1400--1412},
  year={2008},
  doi = {10.2514/1.32836}
}

@article{hou2023optimal,
  title={An optimal geometrical guidance law for impact time and angle control},
  author={Hou, Libing and Luo, Haowen and Shi, Heng and Shin, Hyo-Sang and He, Shaoming},
  journal={IEEE Transactions on Aerospace and Electronic Systems},
  volume={59},
  number={6},
  pages={9821--9830},
  year={2023},
  publisher={IEEE},
  doi = {10.1109/TAES.2023.3305974}
}

@article{hu2018new,
  title={New impact time and angle guidance strategy via virtual target approach},
  author={Hu, Qinglei and Han, Tuo and Xin, Ming},
  journal={Journal of Guidance, Control, and Dynamics},
  volume={41},
  number={8},
  pages={1755--1765},
  year={2018},
  publisher={American Institute of Aeronautics and Astronautics},
  doi = {10.2514/1.G003436}
}

@article{zhang2022virtual,
  title={Virtual target approach-based optimal guidance law with both impact time and terminal angle constraints},
  author={Zhang, Zhihong and Ma, Kemao and Zhang, Gongping and Yan, Liang},
  journal={Nonlinear Dynamics},
  volume={107},
  number={4},
  pages={3521--3541},
  year={2022},
  publisher={Springer},
  doi = {10.1007/s11071-021-07142-3}
}

@article{tekin2017polynomial,
  title={Polynomial shaping of the look angle for impact-time control},
  author={Tekin, Raziye and Erer, Koray S and Holzapfel, Florian},
  journal={Journal of Guidance, Control, and Dynamics},
  volume={40},
  number={10},
  pages={2668--2673},
  year={2017},
  publisher={American Institute of Aeronautics and Astronautics},
  doi = {10.2514/1.G002751}
}

@article{wang2024physics,
  title={A physics-informed indirect method for trajectory optimization},
  author={Wang, Kun and Lu, Fangmin and Chen, Zheng and Li, Jun},
  journal={IEEE Transactions on Aerospace and Electronic Systems},
  volume={60},
  number={6},
  pages={9179--9192},
  year={2024},
  publisher={IEEE},
  doi = {10.1109/TAES.2024.3438687}
}

@article{wu2025nonlinear,
  title={Nonlinear optimal guidance with constraints on impact time and impact angle},
  author={Wu, Fanchen and Chen, Zheng and Shao, Xueming and Wang, Kun},
  journal={Automatica},
  volume={181},
  pages={112500},
  year={2025},
  publisher={Elsevier},
  doi = {10.1016/j.automatica.2025.112500}
}

@article{wang2026ITCG,
  title={Guaranteeing Convergence in Trajectory Shaping Guidance
for Impact Time Control},
  author={Wang, Kun and Ding, Chenxi and Wei, Zhenyu and Wang, Pengyu and Chen, Zheng},
  year={2026},
  volume={191},
  journal = {Automatica},
  pages= {113107},
  doi = {10.1016/j.automatica.2026.113107},
  publisher={Elsevier}
}

\end{document}